\documentclass{gtart}

\def\ifplaintex{\expandafter\ifx\csname documentclass\endcsname\relax}

\def\gtp{{\mathsurround=0pt\it $\cal G\mskip-2mu$eometry \&\ 
$\cal T\!\!$opology $\cal P\!$ublications}}  

\def\recd{{\small Received:\qua\receiveddate\ifx\reviseddate\relax
\else\qquad Revised:\qua\reviseddate\fi\par}} 

\def\lognumber#1{\def\thelognumber{#1}}
\def\volumenumber#1{\def\thevolumenumber{#1}}
\def\volumeyear#1{\def\thevolumeyear{#1}}
\def\papernumber#1{\def\thepapernumber{#1}}
\def\pagenumbers#1#2{\def\startpage{#1}\def\finishpage{#2}}
\def\published#1{\def\publishdate{#1}}

\def\received#1{\def\receiveddate{#1}}

\def\accepted#1{\def\accepteddate{#1}}

\long\def\asciiabstract#1{\long\def\theasciiabstract{#1}}

\let\\\par\let\thelognumber\relax\let\thevolumenumber\relax
\let\thepapernumber\relax\let\thevolumeyear\relax\let\startpage\relax
\let\finishpage\relax\let\publishdate\relax\let\receiveddate\relax
\let\reviseddate\relax\let\accepteddate\relax\let\theasciititle\relax
\let\theasciiauthors\relax
\let\theasciiabstract\relax

\let\theasciiemail\relax

\ifplaintex
\font\logobig=cmssbx10 scaled 3836
\font\logomed=cmssbx10 scaled 2557
\else
\font\logobig=cmssbx10 scaled 4200
\font\logomed=cmssbx10 scaled 2800
\fi

\long\def\makeagttitle{   
\count0=\startpage
\agt\hfill      
\hbox to 45truept{\vbox to 0pt{\vglue -13truept{\logomed A\kern -.37em{\logobig 
T}\kern -.38em G}\vss}\hss}
\break
{\small Volume \thevolumenumber\ (\thevolumeyear)
\startpage--\finishpage\nl
Published: \publishdate}

\vglue .25truein

{\parskip=0pt\leftskip 0pt plus
1fil\def\\{\par\smallskip}{\Large\bf\thetitle}\par\medskip} \vglue
0.05truein

{\parskip=0pt\leftskip 0pt plus 1fil\def\\{\par}{\sc\theauthors}
\par\medskip}%
 
\vglue 0.03truein 

{\small\leftskip 25truept\rightskip 25truept{\bf Abstract}\stdspace\theabstract

{\bf AMS Classification}\stdspace\theprimaryclass
\ifx\thesecondaryclass\relax\else; \thesecondaryclass\fi\par
{\bf Keywords}\stdspace \thekeywords\par}\vglue 7truept

}   

\ifplaintex
\font\phead=cmsl9 scaled 950
\font\pnum=cmbx10 scaled 913
\font\pfoot=cmsl9 scaled 950
\headline{\vbox to 0pt{\vskip -4.5mm\line{\small\phead\ifnum
\count0=\startpage ISSN 1472-2739 (on-line) 1472-2747 (printed)
\hfill {\pnum\folio}\else\ifodd\count0\def\\{ }%
\ifx\theshorttitle\relax\thetitle\else\theshorttitle\fi\hfill{\pnum\folio}
\else\def\\{ and }{\pnum\folio}\hfill\ifx\theshortauthors\relax\theauthors
\else\theshortauthors\fi\fi\fi}\vss}}
\footline{\vbox to 0pt{\vglue 0mm\line{\small\pfoot\ifnum\count0=\startpage
\copyright\ \gtp\hfill\else
\agt, Volume \thevolumenumber\ (\thevolumeyear)\hfill\fi}\vss}}
\else
\font=cmsl9 scaled 1050
\font\lnum=cmbx10 
\font=cmsl9 scaled 1050
\makeatletter
\def\@oddhead{{\small\lhead\ifnum\count0=\startpage ISSN 1472-2739 
(on-line) 1472-2747 (printed)\hfill {\lnum\number\count0}\else\ifodd\count0
\def\\{ }\ifx\theshorttitle\relax \thetitle \else\theshorttitle\fi\hfill
{\lnum\number\count0}\else\def\\{ and }{\lnum\number\count0}
\hfill\ifx\theshortauthors\relax 
\theauthors\else\theshortauthors\fi\fi\fi}}\def\@evenhead{\@oddhead}
\def\@oddfoot{\small\lfoot\ifnum\count0=\startpage\copyright\ \gtp\hfill\else
\agt, Volume \thevolumenumber\ (\thevolumeyear)\hfill\fi}
\def\@evenfoot{\@oddfoot}
\makeatother
\fi
\let\maketitlepage\makeagttitle
\let\makeshorttitle\maketitlepage
\let\maketitle\maketitlepage

\newwrite\gtoutfile
\long\gdef\makeheadfile{  
{\def\\{, }\def\s{ }
\immediate\openout\gtoutfile head.xxx
\immediate\write\gtoutfile{To: math@arxiv.org}
\immediate\write\gtoutfile{Subject: put OR rep NNNNN:ppppp}
\immediate\write\gtoutfile{--text follows this line--}
\immediate\write\gtoutfile{Proxy-for: \ifx\theasciiauthors\relax
\theauthors\else\theasciiauthors\fi\s<\ifx\theasciiemail\relax\theemail\else\theasciiemail\fi>}
\immediate\write\gtoutfile{\noexpand\\}
\immediate\write\gtoutfile{Authors: \ifx\theasciiauthors\relax
\theauthors\else\theasciiauthors\fi}
{\def\\{ }\immediate\write\gtoutfile{Title: \ifx\theasciititle\relax
\thetitle\else\theasciititle\fi}}
\immediate\write\gtoutfile{Subj-class: GT or SG, GR etc}
\immediate\write\gtoutfile{MSC-class: \theprimaryclass\ifx\thesecondaryclass\relax\else, \thesecondaryclass\fi}
\immediate\write\gtoutfile{Journal-ref: Algebr. Geom. Topol. \thevolumenumber\s
(\thevolumeyear) \startpage-\finishpage}
\immediate\write\gtoutfile{Comments: Published by Algebraic and
Geometric Topology at}
\immediate\write\gtoutfile{\s\s\s  http://www.maths.warwick.ac.uk/agt/AGTVol\thevolumenumber/agt-\thevolumenumber-\thepapernumber.abs.html}
\immediate\write\gtoutfile{\noexpand\\}
\immediate\write\gtoutfile{}
\ifx\theasciiabstract\relax
\immediate\write\gtoutfile{\theabstract}\else
\immediate\write\gtoutfile{\theasciiabstract}\fi
\immediate\write\gtoutfile{}
\immediate\write\gtoutfile{\noexpand\\}
\immediate\write\gtoutfile{}
\immediate\closeout\gtoutfile}}  

\def\maketitlepage{\makeagttitle\makeheadfile}
\let\makeshorttitle\maketitlepage
\let\maketitle\maketitlepage

\lognumber{6}
\volumenumber{3}
\volumeyear{2003}
\papernumber{6}
\published{20 February 2003}
\pagenumbers{147}{154}
\received{27 January 2003}
\accepted{11 February 2003}

\usepackage{amsmath, amssymb, graphicx}

\newtheorem{thm}{Theorem}[section]

\theoremstyle{definition}
\newtheorem{defn}[thm]{Definition}

\newtheorem{rem}[thm]{Remark}

\newcommand{\field}[1]{\mathbb{#1}}
\newcommand{\integers}{\ensuremath{\field{Z}}}
\newcommand{\naturals}{\ensuremath{\field{N}}}
\newcommand{\reals}{\ensuremath{\field{R}}}
\newcommand{\Euclidean}{\ensuremath{\field{E}}}

\begin{document}

\title
[A flat that is not the limit of periodic flats]
{A flat plane that is not the limit of periodic\\flat planes}
\author{Daniel T. Wise}
\address{Department of Mathematics and Statistics, McGill 
University\\Montreal, Quebec H3A 2K6, Canada}
\email{wise@math.mcgill.ca}


\date{\today}

\begin{abstract}
We construct a compact nonpositively curved squared $2$-com\-p\-lex
whose universal cover contains a flat plane that is not the limit
of periodic flat planes.
\end{abstract}

\asciiabstract{We construct a compact nonpositively curved squared
2-complex whose universal cover contains a flat plane that is not the
limit of periodic flat planes.}

\maketitle

\section{Introduction}

Gromov raised the question of which ``semi-hyperbolic spaces''
have the property that their flats can be approximated by periodic flats
\cite[\S6.B$_3$]{Gromov93}.
In this note we construct an example of a
compact nonpositively curved squared $2$-complex $Z$
whose universal cover $\tilde Z$ contains an isometrically embedded
flat plane that is not the limit of a sequence of periodic flat planes.

A flat plane $\Euclidean \hookrightarrow \tilde Z$ is {\em periodic}
if the map $\Euclidean\looparrowright Z$ factors as
$\Euclidean\rightarrow T\rightarrow Z$ where $\Euclidean \rightarrow T$
is a covering map of a torus $T$.
Equivalently, $\pi_1Z$ contains a subgroup isomorphic to
 $\integers\times\integers$
which stabilizes $\Euclidean$ and acts cocompactly on it.
A flat plane $f\co\Euclidean\hookrightarrow \tilde Z$ is the
{\em limit of periodic flat planes}
if there is a sequence of periodic flat planes
$f_i\co\Euclidean\hookrightarrow  \tilde Z$ which converge
pointwise to $f\co\Euclidean\rightarrow \tilde Z$.
In our setting, $\tilde Z$ is a $2$-dimensional complex,
and so $\Euclidean\hookrightarrow \tilde Z$ is the limit of periodic flat planes
 if and only if every compact subcomplex of $\Euclidean$ is contained in a periodic flat plane.

In Section~\ref{sec:X antitorus} we describe a compact
nonpositively curved $2$-complex $X$ whose universal cover
contains a certain aperiodic plane called an ``anti-torus''.
In Section~\ref{sec:Z} we construct $Z$ from $X$ by strategically
gluing tori and cylinders to $X$ so that $\tilde Z$ contains
a flat plane which is a mixture of the anti-torus and periodic planes.
This flat plane is not approximable by periodic flats because
it contains a square that does not lie in any periodic flat.
Our example $Z$ is a $K(\pi,1)$ for a negatively
curved square of groups, and in Section~\ref{sec:polygon}
we describe an interesting related triangle of groups.

\section{The anti-torus in $X$}\label{sec:X antitorus}
\subsection{The $2$-complex $X$}
\begin{figure}\centering
\includegraphics[width=\textwidth]{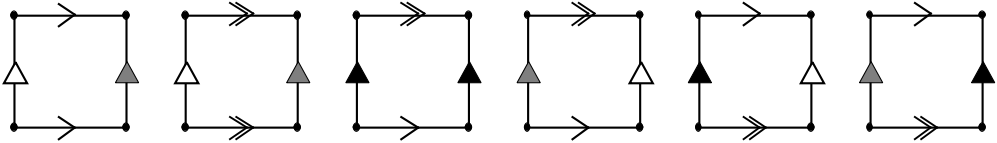}
\caption{The figure above indicates the gluing pattern for the six squares of
$X$. The three vertical edges colored white, grey, and  black
are denoted $a$, $b$, and  $c$ respectively. The two horizontal
edges, single and double arrow, are denoted $x$ and $y$
respectively.\label{fig:sixsquares}}
\end{figure}

Let $X$ denote the complex consisting of the six squares
indicated in Figure~\ref{fig:sixsquares}. The squares are glued together as
indicated by the oriented labels on the edges. Note that $X$ has a unique
0-cell, and that the notion of vertical and horizontal
are preserved by the edge identifications.
Let $H$ denote the subcomplex consisting of the
$2$~horizontal edges, and let  $V$ denote  the subcomplex
consisting of the $3$~vertical edges.

The complex $X$, which was first studied in \cite{Wise96Thesis},
has a number of interesting properties that we record here:
The link  of the unique $0$-cell in $X$ is a complete bipartite graph.
It follows that the universal cover $\widetilde X$ is the product of two trees
$\tilde H\times \tilde V$ where $\tilde H$ and $\tilde V$ are the universal covers
of $H$ and $V$.
In particular, the link contains no cycle of length~$<4$ and so
$X$ satisfies the combinatorial nonpositive curvature condition
for squared $2$-complexes \cite{Gromov87,BridsonHaefliger} which is a special
case of the $C(4)$-$T(4)$ small-cancellation condition \cite{LS77}.

The $2$-complex $X$ was used in
\cite{Wise96Thesis} to produce the first examples of non-residually finite
groups which are fundamental groups of spaces with the above properties.
The connection to finite index subgroups arises because
while $\widetilde X$ is isomorphic to the cartesian product of two trees,
 $X$ does not have a finite cover which is the
product of two graphs.

\subsection{The anti-torus $\Pi$}

\begin{figure}\centering
\includegraphics[width=.4\textwidth]{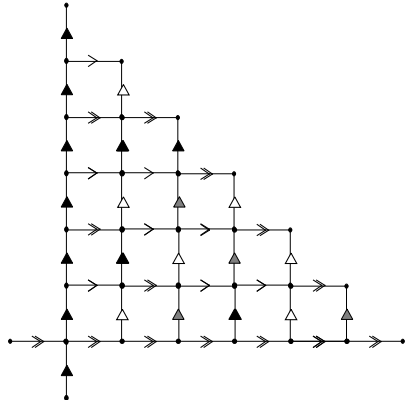}
\caption{ The Anti-Torus $\Pi$: The plane $\Pi$ above is the
convex hull of two periodically labeled lines in $\tilde X$. A
small region of the northeast quadrant has been tiled by the
squares of $X$. \label{fig:plane}}
\end{figure}

The exotic behavior of $X$ can be attributed to
the existence of a strangely aperiodic plane $\Pi$ in $\widetilde X$
that we shall now describe.
Let $\tilde x\in \widetilde X^0$ be the basepoint of $\widetilde X$.
Let $c^\infty$ denote
the infinite periodic vertical line in $\tilde X$ which is the
based component of the preimage of the loop labeled by $c$ in $X$.
Define $y^\infty$ analogously.
Let $\Pi$ denote the convex hull
in $\tilde X$ of the infinite geodesics labeled by $c^\infty$ and
$y^\infty$, so $\Pi=y^\infty\times c^\infty$. The plane  $\Pi$ is tiled by the six orbits of
squares in $\tilde X$ as in Figure~\ref{fig:plane}.
The reader can extend $c^\infty\cup y^\infty$
to a flat plane by successively adding squares wherever
there is a pair of vertical and horizontal edges meeting at a vertex.
From a combinatorial point of view,
the existence and uniqueness of this extension is guaranteed by
the fact that the link of $X$ is a complete bipartite graph.

The ``axes'' $c^\infty$ and $y^\infty$ of $\Pi$ are obviously periodic,
and using that $X$ is compact, it is easy to verify that
for any $n\in\naturals$, the infinite strips
$[-n,n]\times \reals$ and $\reals\times [-n,n]$ are periodic.
However, the period of these infinite strips increases exponentially
with~$n$.
Thus, the entire plane $\Pi$ is aperiodic.
Note that to say that $[-n,n]\times \reals$ is {\em periodic}
means that the immersion $\big([-n,n]\times \reals\big)\looparrowright X$
factors as $\big([-n,n]\times \reals\big)\rightarrow C \looparrowright X$
where $\big([-n,n]\times \reals\big)\rightarrow C$ is the universal
covering map
of a cylinder. The map $\Pi\looparrowright X$ is {\em aperiodic} in the sense
that it does not factor through an immersed torus.

We conclude this section by giving a brief explanation
of the aperiodicity of $\Pi$.
A complete proof that $\Pi$ is aperiodic is given in \cite{Wise96Thesis}.
Let $W_n(m)$ denote the word corresponding to the length~$n$ horizontal
 positive path in $\Pi$  beginning at the endpoint of the vertical path $c^m$.
 Thus, $W_n(m)$ is the label of the side opposite $y^n$
in the rectangle which is the combinatorial convex hull of $y^n$ and $c^m$.
Equivalently, $W_n(m)$ occupies the interval $\{m\}\times [0,n]$.
For each $n$, the words
$\{W_n(m) \mid  0 \leq m \leq 2^n -1 \}$ are all distinct!
Consequently every  positive length~$n$ word  in $x$ and $y$ is $W_n(m)$ for some $m$.
This implies that the infinite strip $[0,n]\times \reals$
has period $2^n$, and in particular $\Pi$ cannot be periodic.

We refer to $\Pi$ as an {\em anti-torus} because the aperiodicity of
$\Pi$ implies that $c$ and $y$ do not have nonzero powers which commute.
Indeed, if $c^p$ and $y^q$ commuted for $p,q\neq0$ then the flat torus
theorem (see \cite{BridsonHaefliger}) would imply that
$c^\infty$ and $y^\infty$ meet in a periodic flat plane,
which would contradict that $\Pi$ is aperiodic.

 \section{The $2$-complex $Z$ with a nonapproximable flat}
 \label{sec:Z}
\begin{figure}\centering
\includegraphics[width=.2\textwidth]{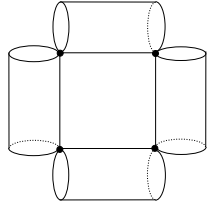}
\caption{ The complex $Y$ is formed by gluing four cylinders to a square.
\label{fig:Y}}
\end{figure}

 We first construct a new complex $Y$ as follows: Start with a square
$s$, and then attach four cylinders each of which is isomorphic
to $S^1\times I$. One such cylinder is attached  along each side of $s$.
The resulting complex $Y$ containing exactly five squares is illustrated in
Figure~\ref{fig:Y}.

Let $T^2$ denote the torus $S^1\times S^1$ with the usual product cell structure
consisting of one $0$-cell, two $1$-cells, and a single square $2$-cell.
We let $\widetilde T^2$ denote the universal cover and we shall identify
$\widetilde T^2$ with $\reals^2$.

\begin{figure}\centering
\includegraphics[width=.4\textwidth]{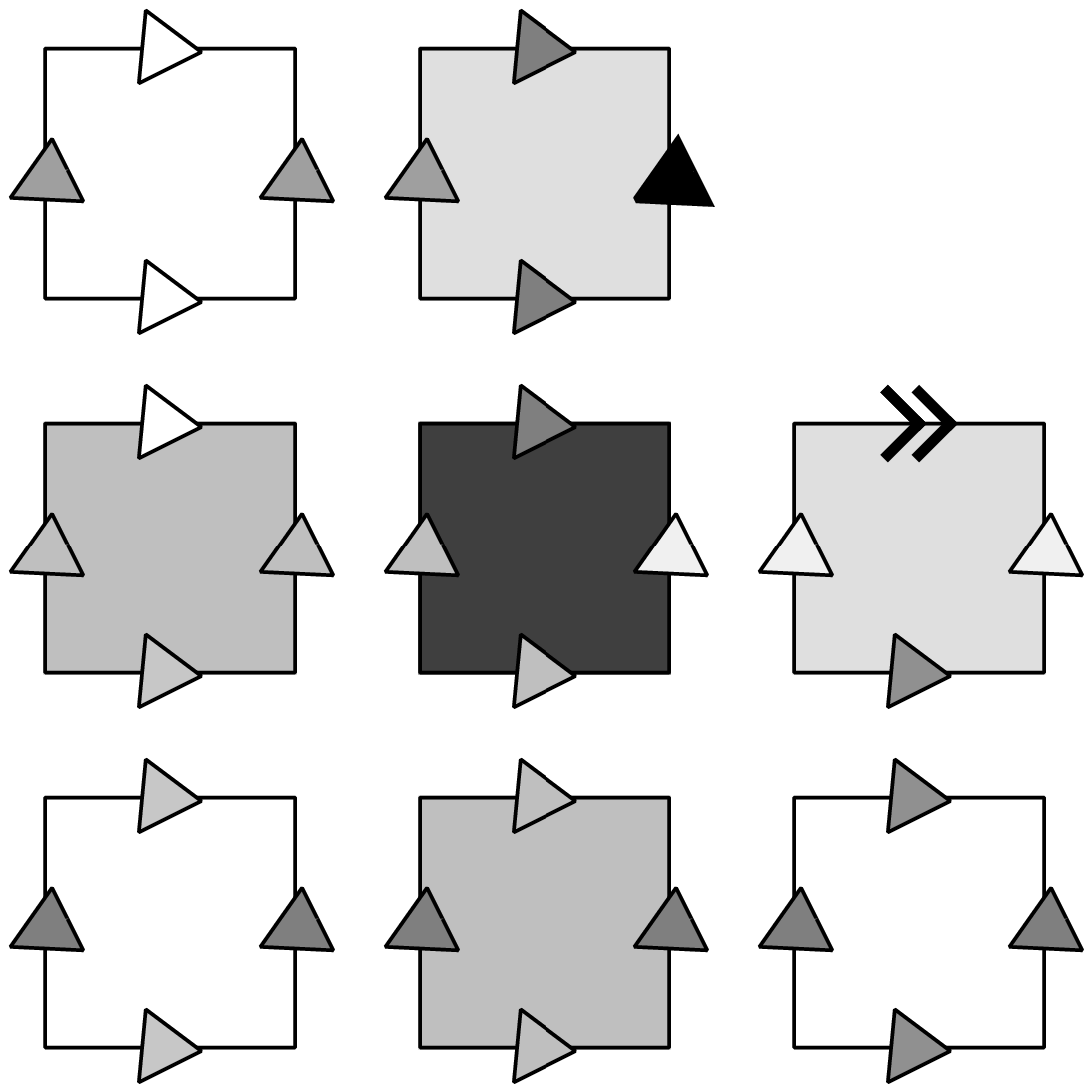}
\hspace{5mm}
\includegraphics[width=.5\textwidth]{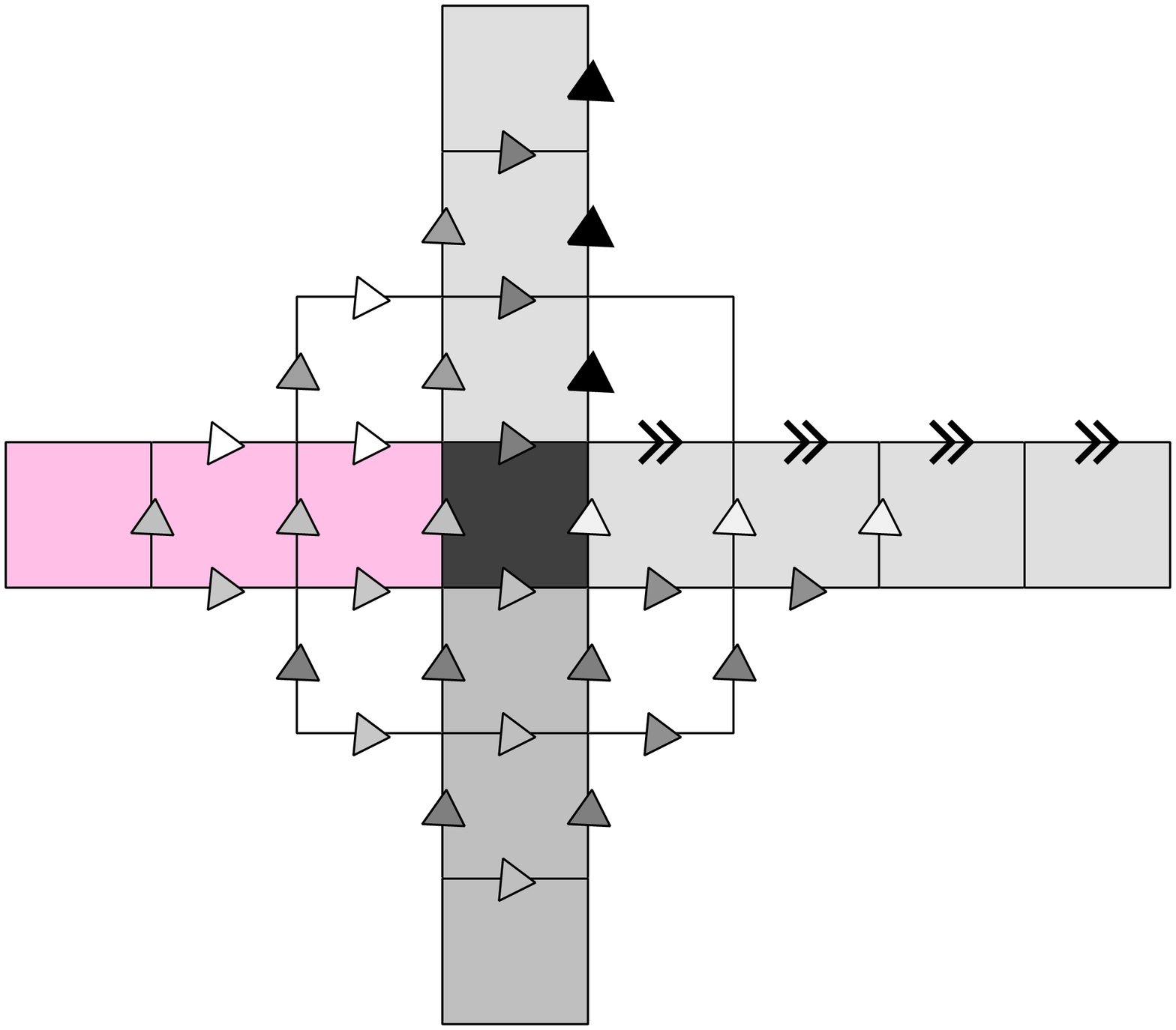}
\caption{
\label{fig:compandplane} $Z-X$ and $Z$: The eight squares of the
figure on the left are glued together following the gluing pattern
to form $Z-X$. To form $Z$, we add a copy of $X$ at the $NE$
corner, identifying the loops in $X$ labeled by $c$ and $y$, with
the black single and double arrows of the diagram. The figure on
the right represents an infinite cross whose convex hull in $Z$
is not approximable by any periodic plane. Note that while the NW,
SW, and SE quarters of this plane are periodic,  the NE quarter
is an aperiodic quarter of $\Pi$.}
\end{figure}

At each corner of $s \subset Y$, there is a pair of intersecting circles in $Y^1$,
which are boundary circles of distinct cylinders.
Note that they meet at an angle of $\frac{3\pi}2$ in $Y$.
At each of three (NW, SW, \& SE)
 corners of $s \subset Y$ we attach a copy of $T^2$
 by identifying the pair of circles in the  $1$-skeleton of $T^2$
 with the pair of intersecting circles noted above at
the respective corner of $s$.
At the fourth (NE) corner of $s$, we attach a copy of the complex
$X$. Here we identify the pair of circles meeting at the corner of
$s$ with the pair of perpendicular circles $c$ and $y$ of $X$.
We denote the resulting complex by $Z$. Thus, $Z =   T^2 \cup T^2 \cup
T^2 \cup Y \cup X$. See Figure~\ref{fig:compandplane} for a
depiction of the $8$ squares of $Z - X$ and their gluing patterns.

\begin{defn}[Infinite cross]
An {\em infinite cross}
is a squared $2$-complex isomorphic to the subcomplex of $\widetilde T^2$
consisting of $\big([0,1] \times \reals\big)
\cup \big(\reals \times [0,1]\big)$.
The {\em base square} of the infinite cross is the square $[0,1]\times[0,1]$.
\end{defn}

\begin{figure}\centering
\includegraphics[width=.5\textwidth]{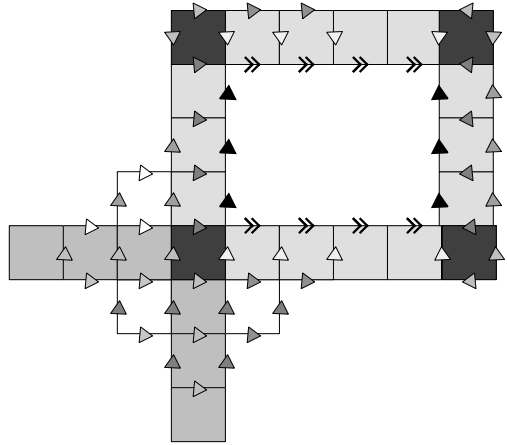}
\caption{
The figure above
illustrates one of the four possible contradictions which explain
why no  periodic plane contains the square $s$.
\label{fig:contradiction}}
\end{figure}

{\bf  The planes containing $s$:}\qua
Observe that $Y$ contains various immersions of an infinite
cross whose base square maps to $s$.
In particular, there are exactly 16~distinct immersed infinite crosses
$C\looparrowright Y$ that pass through $s$ exactly once.
Each of these infinite crosses extends uniquely to an immersed flat plane
in $Z$.
Each such flat plane fails to be periodic because its four quarters map to distinct
parts of $Z$.
Our main result is that these immersed flat planes
are not approximable by periodic flat planes because of the following:

\begin{thm}[No periodic approximation] There is no
immersion of a torus $T^2 \rightarrow Z$ which contains $s$.
Equivalently, there is no periodic plane in $\tilde Z$ containing
$\tilde s$.
\end{thm}
\begin{proof}
 We argue by contradiction. Suppose that there is an immersed periodic plane
$\Omega$ containing $s$.
We shall now produce a rectangle as in Figure~\ref{fig:contradiction}
that will yield a contradiction.
We may assume that a copy of $s$ in $\Omega$ is oriented
as in Figure~\ref{fig:compandplane}.
We begin at this copy of $s$ and travel
north inside the northern cylinder until we reach another
copy $s_n$ of $s$. The existence of $s_n$ is guaranteed by our assumption
that $\Omega$ is periodic.
Similarly, we travel east from $s$ to reach a square $s_e$.
Travelling north from $s_e$ and east from $s_n$, we trace out the
boundary of a rectangle whose NE corner is a  square $s_{ne}$
(see Figure~\ref{fig:contradiction}).

This yields a contradiction because the inside of this
rectangle is tiled by squares in $X$, yet the boundary of this
rectangle is a commutator $\big [ c^{\pm n}, y^{\pm m} \big] $.
As explained in Section~\ref{sec:X antitorus},
such a word cannot be trivial in $\pi_1X$
 because of the anti-torus.
\end{proof}

\begin{rem}
Using an argument similar to the above proof, one can show
that these sixteen planes are the {\it only} flat planes in $\tilde Z$ containing
$\tilde s$.
One considers the pair of ``axes''
intersecting at $\tilde s$ in a plane containing $\tilde s$. If this plane is
different from each of the $16$ mentioned above, then some translate
of $\tilde s$ must appear along one of these ``axes''. The infinite strip
in the plane whose corners are these two $s$ squares
yields a contradiction similar to the one obtained above.
\end{rem}

\begin{rem}
While $X$ is a rather pathological complex,
we note that every flat plane in $\tilde X$ is the limit of periodic flat planes.
Indeed this holds for any compact $2$-complex $X$ whose universal cover
is isomorphic to the product of two trees
\cite{Wise96Thesis}.
\end{rem}

\section{Polygons of groups}
\label{sec:polygon}
\subsection{The algebraic angle versus the geometric angle}
Since the elements $c$ and $y$ have axes which
intersect perpendicularly in a plane in $\tilde X$, the
natural geometric angle between the subgroups $\langle c\rangle$ and
$\langle y\rangle$ of $\pi_1X$ is $\frac\pi2$.
However, the algebraic Gersten-Stallings angle (see \cite{Stallings91}) between
these subgroups is $\leq \frac\pi3$.
To see this, we must show that there is no non-trivial
relation of the form $c^k y^l c^m y^n = 1$.

Since $\tilde X$ is isomorphic to the cartesian product $\tilde V\times \tilde H$,
of two trees and $c$ and $y$ correspond to distinct factors,
it follows that the only relations that must be checked are
rectangular
 (i.e., $\vert k\vert = \vert m\vert$ and
        $\vert l\vert = \vert n\vert$ ).
        However, these are
easily ruled out by the anti-torus $\Pi$ and the fact that $X$ is nonpositively curved.

\subsection{Square of groups and triangle of groups}
The complex $Z$ can be  thought of in a natural way
as a $K(\pi,1)$ for a negatively curved {\it square of groups}
(see \cite{Stallings91,Haefliger91,Corson92}) with cyclic edge groups
and trivial face group.

Because  the algebraic angle between $\langle c \rangle$ and
$\langle y \rangle$ in $\pi_1X$ is $\leq \frac\pi3$, it is tempting to form an analogous
nonpositively curved triangle of groups $D$.
The face group of $D$ is trivial, the edge groups
of $D$ are cyclic, the vertex groups of $D$ are isomorphic to $\pi_1X$, and
each edge group of $D$ is embedded on one (clockwise) side as
$\langle c \rangle$ and on the other (counter-clockwise) side as
$\langle y \rangle$. This can be done so that the resulting
triangle of groups $D$ has $\integers_3$ symmetry.
The tension between the algebraic and geometric angles
should endow $\pi_1D$ with some interesting properties.
For instance, I suspect that $\pi_1D$ fails to be the fundamental group
of a compact nonpositively curved space, but it fails for reasons different
from the usual types of problems.

{\bf Acknowledgments}\qua I am grateful to the referee for the helpful
corrections.  This research was supported by an NSERC grant.

\Addresses\recd
\end{document}